\documentclass[preprint,12pt]{elsarticle}

\usepackage{amssymb,amsthm}

\newtheorem{Def}{Definition}[section]
\newtheorem{Th}{Theorem}[section]
\newtheorem{Prop}{Proposition}[section]
\newtheorem{Not}{Notation}[section]
\newtheorem{Lemma}{Lemma}[section]
\newtheorem{lemma}{Lemma}
\newtheorem{Rem}{Remark}[section]
\newtheorem{Cor}{Corollary}[section]
\def\s{\section}                        \def\ss{\subsection}
\def\d{\begin{Def}}                     \def\t{\begin{Th}}
\def\p{\begin{Prop}}                    \def\n{\begin{Not}}
\def\la{\begin{Lemma}}                  \def\r{\begin{Rem}}
\def\c{\begin{Cor}}                     \def\ee{\begin{equation}}
\def\aa{\begin{eqnarray}}               \def\y{\begin{eqnarray*}}
\def\bd{\begin{description}}

\def\ed{\end{Def}}                      \def\et{\end{Th}}
\def\ep{\end{Prop}}                     \def\en{\end{Not}}
\def\el{\end{Lemma}}                    \def\er{\end{Rem}}
\def\ec{\end{Cor}}                      \def\eee{\end{equation}}
\def\eaa{\end{eqnarray}}                \def\ey{\end{eqnarray*}}
\def\ebd{\end{description}}

\def\nn{\nonumber}

                        \def\lb{\label}

                         \def\R{{\bf R}}
                         
                         \def\Ua{\mathcal{U}}

\def\g{{\bf g}}

                         \def\Sa{{\rm S}}

                   \def\tr{{\rm tr}}

\begin{document}

\begin{frontmatter}

%% Title, authors and addresses

%% use the tnoteref command within \title for footnotes;
%% use the tnotetext command for the associated footnote;
%% use the fnref command within \author or \address for footnotes;
%% use the fntext command for the associated footnote;
%% use the corref command within \author for corresponding author footnotes;
%% use the cortext command for the associated footnote;
%% use the ead command for the email address,
%% and the form \ead[url] for the home page:
%%
%% \title{Title\tnoteref{label1}}
%% \tnotetext[label1]{}
%% \author{Name\corref{cor1}\fnref{label2}}
%% \ead{email address}
%% \ead[url]{home page}
%% \fntext[label2]{}
%% \cortext[cor1]{}
%% \address{Address\fnref{label3}}
%% \fntext[label3]{}

\title{Gradient estimates for $u_t=\Delta F(u)$ on manifolds and some Liouville-type theorems}

%% use optional labels to link authors explicitly to addresses:
%% \author[label1,label2]{<author name>}
%% \address[label1]{<address>}
%% \address[label2]{<address>}

\author{Xiangjin Xu \fnref{fn1}}
\ead{xxu@math.binghamton.edu}

\address{Department of Mathematical Sciences, Binghamton University - SUNY\\ Binghamton, NY 13902, U.S.A}

%\cortext[cor1]{Corresponding author}
\fntext[fn1]{The author is partially supported by NSF grants NSF-DMS-0602151 and NSF-DMS-0852507, and Harpur College Grant in Support of Research, Scholarship and Creative Work from Harpur College, Binghamton University. }

\begin{abstract} In this paper, we first prove a localized Hamilton-type gradient estimate for the positive solutions of Porous Media type equations:
$$u_t=\Delta F(u),$$
with $F'(u) > 0$, on a complete Riemannian manifold with Ricci curvature bounded from below. In the second part, we study Fast Diffusion Equation (FDE) and Porous Media Equation (PME):
$$u_t=\Delta (u^p),\qquad p>0,$$
and obtain localized Hamilton-type gradient estimates for FDE and PME in a larger range of $p$ than that for Aronson-B\'enilan estimate, Harnack inequalities and Cauchy problems in the literature. Applying the localized gradient estimates  for FDE and PME, we prove some Liouville-type theorems for positive global solutions of FDE and PME on noncompact complete manifolds with nonnegative Ricci curvature, generalizing Yau$¡¯$s celebrated Liouville theorem for positive harmonic functions.

%% Text of abstract

\end{abstract}

\begin{keyword}
Localized Hamilton-type gradient estimate, Liouville-type Theorem, Fast Diffusion Equation, Porous Media Equation

%% keywords here, in the form: keyword \sep keyword
%% MSC codes here, in the form: \MSC code \sep code
%% or \MSC[2008] code \sep code (2000 is the default)
\MSC[2010] 53C21\sep 35K05

\end{keyword}

\end{frontmatter}

%%
%% Start line numbering here if you want
%%
% \linenumbers

%% main text
%\section{}
%\label{}

\s{Introduction}

%\vspace{2mm}

The goal of the paper is to establish a localized Hamilton-type gradient estimate for the positive solutions of Porous Media type equations, which are degenerate parabolic equations in general,
\begin{eqnarray}
u_t=\Delta F(u)\lb{Eq1}
\end{eqnarray}
on a complete Riemannian manifold $(M^n, \g)$ with $Ric(M^n) \geq -k$ for some $k \geq 0$. Here $F\in C^2(0,\infty)$, $F'> 0$, and $\Delta$ is the Laplace-Beltrami operator of the metric $\g$. The equation (\ref{Eq1}) is a nonlinear version of the classical heat equation (case $F(u)=u$). Porous Media Equation (PME for short) (case $F(u)=u^p,\; p>1$) has arisen in different applications to model diffusive phenomena like
groundwater infiltration (Boussinesq¡¯s model, $p = 2$), flow of gas in porous media (Leibenzon-Muskat model, $p \geq 2$), heat radiation in plasmas ($m > 4$), and others. The mathematical theory started in the 1950¡¯s and got momentum in recent decades as a nonlinear diffusion problem with interesting geometrical aspects (free boundaries) and peculiar functional analysis. We refer to the monograph \cite{V2} for an account of the rather complete theory concerning existence, uniqueness, regularity and asymptotic behavior of PME. Some of the existence, uniqueness and regularity properties hold true for the so-called Fast Diffusion Equation (FDE for short)(case $F(u)=u^p,\; p\in (0,1)$). FDE appears in plasma physics and in geometric flows such as the Ricci flow on surfaces and the Yamabe flow.

It is well known that, in the study of geometric analysis and other elliptic or parabolic equations, the gradient estimate and the Harnack inequality play an important role.  The Li-Yau estimate and the Harnack inequality in the fundamental paper \cite{LY}, where Li and Yau studied the heat equation on general Riemannian manifolds, have a tremendous impact on the field of geometric analysis. Since 1970s, Aronson-B\'enilan estimate and the Harnack-type inequalities have been widely studied for PME and FDE defined on the whole Euclidean space, cf. \cite{Ar}, \cite{ArB}, \cite{AB}, \cite{HP}. Recently, Lu, Ni, V\'azquez and Villani \cite{LNVV} studied PME and FDE on manifolds and got some localized Aronson-B\'enilan estimates. The first Harnack-type inequality dealing with the Porous Media type  equation (\ref{Eq1}) was attributed to S.T. Yau \cite{Y}.

{\bf Theorem A (S.T. Yau \cite{Y}). } {\it Let $M^n$ be a compact Riemannian manifold without boundary, $Ricci(M) \geq 0$. Suppose that $F \in C^2(0,\infty)$ with $F' > 0$, $c(t) \in C^1(0,\infty)$, and $u$ is any positive solution of the degenerate parabolic equation
\ee
u_t=\Delta F(u)\nn
\eee
on $M^n$. Let $\alpha \neq 0$ be an arbitrary constant. Define a function $G$ on $(0,\infty)$ by $G'(s) = F'(s)/s$, and we abbreviate $G = G(u)$, $F^{\kappa} = F^{\kappa}(u)$, $\kappa = 0, 1, 2$.

If the conditions below are satisfied:

(A). $|\nabla G|^2-\alpha G_t-c(t)\leq 0$ at $t= 0$;

(B). (nonlinear condition) the following quadratic inequality holds true for all $x \geq 0$
\begin{eqnarray*}
0&\geq& \frac{1-\alpha}{\alpha^2}\Bigg(\alpha u F''- \frac{2(1-\alpha)}{n}F'\Bigg)x^2  +\Bigg(\frac{4(1-\alpha)}{n\alpha^2} -\frac{u F''}{F'}\Bigg)c(t)x \\
&& -\Bigg(\frac{2}{n}+\alpha\frac{u F''}{F'}\Bigg)\frac{c^2(t)}{\alpha^2F'} - c'(t)
\end{eqnarray*}
then we have for all $t > 0$ that
\begin{eqnarray*}
|\nabla G|^2-\alpha G_t-c(t)\leq 0.
\end{eqnarray*}
}

For the heat equation on compact manifolds without boundary, Hamilton \cite{H} studied another type of gradient estimates as:

%\vspace{2mm}

{\bf Theorem B (Hamilton \cite{H}).} {\it Let $M^n$ be a compact manifold without boundary and with $Ric(M^n) \geq -k$ for some $k \geq 0$. Let $u$ be a smooth positive solution of the heat equation with $u \leq M$ for all $(x, t) \in M^n\times (0,\infty)$. Then
\begin{eqnarray*}
\frac{|\nabla u|^2}{u^2} \leq \Bigg(\frac{1}{t}+2k\Bigg)\ln\frac{M}{u}.
\end{eqnarray*}}

%\vspace{2mm}

The Hamilton-type gradient estimate takes up a significant position in the study of the heat equation. However, the classical Hamilton$'$s estimate is a global result which requires the heat equation to be posed on compact manifold without boundary. Recently, a localized Hamilton type gradient estimate was proven by  Souplet and Zhang \cite{SZ}, which can be viewed as a combination of Li-Yau estimate and Hamilton$'s$ gradient estimate.

{\bf Theorem C (Souplet and Zhang \cite{SZ}):} {\it Let $M^n$ be a complete Riemannian manifold with dimension $n\geq 1$, $Ric(M^n)\geq -k$, $k \geq 0$. Suppose that $F\in C^2(0,\infty)$ with $F' > 0$, and $u$ is any positive solution of the degenerate parabolic equation (\ref{Eq1}) in $Q_{R,T} \equiv B(x_0,R)\times [t_0-T , t_0]\subset M^n\times(-\infty,\infty)$. Suppose also that $u\leq M$ in $Q_{R,T}$. Then there exists a dimensional constant $C$ such that
\begin{eqnarray*}
\sup_{(x, t)\in Q_{R/2,T /2}}\frac{|\nabla_x u(x,t)|}{u(x,t)}\leq C\Bigg( \frac{1}{R} +\frac{1}{\sqrt{T}}+\sqrt{k} \Bigg)\Bigg(1+\ln\frac{M}{u}\Bigg).
\end{eqnarray*}

Moreover, if $M^n$ has nonnegative Ricci curvature and $u$ is any positive solution of the heat equation on $M \times (0, \infty)$, then there exist dimensional constants $C_1, C_2$ such that
\begin{eqnarray*}
\frac{|\nabla_x u(x,t)|}{u(x,t)}\leq C_1\frac{1}{\sqrt{t}}\Bigg(C_2+\ln\frac{u(x,2t)}{u(x,t)}\Bigg).
\end{eqnarray*}
for all $x \in M^n$ and $t > 0$.}\\

It is natural to seek a localized Hamilton-type gradient estimate for Porous Media type equation (\ref{Eq1}) as Souplet and Zhang \cite{SZ} did for the heat equation on a complete manifold. Recently, Ma,  Zhao and Song \cite{MZS} proved  a localized Hamilton-type gradient estimate for the equation (\ref{Eq1}) under some strong assumptions (Theorem  7 in \cite{MZS}), where the gradient estimate for FDE ($u_t=\Delta (u^p)$) holds only for   dimension $n=2,3$ with $p\in(1-1/\sqrt{n},1)$ (Corollary 9 in \cite{MZS}). Such restrictions on $p$ and $n$ are unnatural and inadequate for applications since the mathematical theory of PME and FDE based on a priori estimates such as Aronson-B\'enilan estimate and others (cf. \cite{Ar}, \cite{AB}, \cite{HP}, \cite{LNVV}, etc), applies to all positive smooth solutions of PME and FDE on the condition that $p > p_c := 1 - 2/n$ for any dimension $n$.

The first main result of this paper is the following localized Hamilton-type gradient estimate for the Porous Media type equation (\ref{Eq1}), which generalizes {\bf Theorem C} of Souplet and Zhang \cite{SZ} for the heat equation, under suitable conditions on $F$, as {\bf Theorem A} on the Harnack-type inequality for the Porous Media type  equation (\ref{Eq1}) by Yau \cite{Y}:

\t\label{Thm1.1} {\bf (Gradient Estimate).}  Let $M^n$ be a complete Riemannian manifold with dimension $n\geq 1$, $Ric(M^n)\geq -k$, $k \geq 0$. Suppose that $F\in C^2(\R^+)$ with $F' > 0$, and $u$ is any positive solution of the equation (\ref{Eq1}) in $Q_{R,T} \equiv B(x_0,R)\times [t_0-T , t_0]\subset M^n\times(-\infty,\infty)$. Define a function $G$ on $(0,\infty)$ by $G'(s) =F'(s)/s$. Denote $\Ua=range_{(x,t)\in Q_{R,T}} u(x,t) \subset (0,\infty)$. Choose nonnegative constants $K$, $\alpha$, $\delta$ and $\tau$ such that
\begin{eqnarray*}
 K \geq \sup_{s\in\Ua }F'(s), \quad \alpha-\sup_{s\in\Ua }G(s)\geq \delta>0; \quad \tau \geq \sup_{s\in\Ua }|sF''(s)|/F'(s).
\end{eqnarray*}
If there exists a nonnegative constant $\gamma$, such the condition below is satisfied:
\begin{eqnarray*}
(C):\quad 2+\frac{sF''(s)}{F'(s)}\Big(2-(n-1)\frac{sF''(s)}{F'(s)}\times\frac{\alpha-G(s)}{F'(s)}\Big) \geq \gamma>0, \quad \forall\; s\in\Ua,
\end{eqnarray*}
then there exists a constant $C(n, \delta, K, \tau,\gamma)$ depending only on $n$, $\delta$, $K$, $\tau $ and $\gamma$  such that
\begin{eqnarray}
\sup_{(x, t) \in Q_{R/2,T /2}}\frac{|\nabla_x G(u(x,t))|}{\alpha-G(u(x,t))}\leq C(n, \delta, K, \tau,\gamma)\Bigg( \frac{1}{R} +\frac{1}{\sqrt{T}}+\sqrt{k} \Bigg). \lb{Eq2}
\end{eqnarray}
When $n=1$, the Ricci curvature lower bound $k$ vanishes.\\
\et

In the second part of this paper, we study Fast Diffusion Equation (FDE for short) and Porous Media Equation (PME for short):
\ee
u_t=\Delta (u^p), \qquad p>0\lb{Eq3}
\eee
on a complete Riemannian manifold $(M^n, \g)$ with $Ric(M^n) \geq -k$ for some $k \geq 0$.

Firstly  for FDE, i.e. the equation (\ref{Eq3}) with $p<1$, we obtain the localized Hamilton-type gradient estimate, which generalizes the gradient estimate  for FDE (\ref{Eq3}) in \cite{MZS} (Corollary 9 in \cite{MZS}), where the estimate holds only for dimension $n=2\; or \;3$ with  $p\in(1-1/\sqrt{n},1)$,

%\vspace{2mm}

\t\label{Thm1.2}  Let $M^n$ be a complete Riemannian manifold with dimension $n\geq 1$, $Ric(M^n)\geq -k$, $k \geq 0$. Suppose that $u\leq M$ is a positive solution of FDE (\ref{Eq3}) with
\begin{eqnarray*}
1-\frac{4}{n+3}< p< 1,  & for\; n\geq 1,\nn
\end{eqnarray*}
in $Q_{R,T} \equiv B(x_0,R)\times [t_0-T , t_0]\subset M^n\times(-\infty,\infty)$, then there exists a constant $C(n,p)$ depending only on n and p such that

\begin{eqnarray}
\sup_{(x, t) \in Q_{R/2,T /2}}\frac{|\nabla u(x,t)|}{u(x,t)}\leq  C(n,p)\Big(\frac{1}{R}+\frac{M^{(1-p)/2}}{\sqrt{T}}+\sqrt{k} \Big). \lb{Eq4}
\end{eqnarray}
When $n=1$, the Ricci curvature lower bound $k$ vanishes.

\et

%\vspace{2mm}

\r
One should notice that the range of $p$ here is $(1-\frac{4}{n+3}, 1)$, while previous results on Aronson-B\'enilan estimate and the Harnack-type differential inequalities for FDE  (cf. \cite{Ar}, \cite{AB}, \cite{HP}, \cite{LNVV}, etc) require $p\in ((1-\frac{2}{n})_+, 1)$. We can see that for $n\geq 3$, the range of $p$ for our gradient estimate is larger than that in previous results \cite{Ar}, \cite{AB}, \cite{HP}, \cite{MZS}. Our gradient estimate will be a useful tool to study related problem for FDE in this large range of $p$, in which one couldn't deal with in the literature even on $\R^n$.
\er

%\vspace{2mm}

On a complete noncompact manifold with nonnegative Ricci curvature, an immediate application of Theorem \ref{Thm1.2} is the following time-dependent Liouville theorem for positive global solutions of FDE, generalizing Yau¡¯s celebrated Liouville theorem for positive harmonic functions, which states that any positive harmonic function on a noncompact manifold with nonnegative Ricci curvature is a constant function.

\t\lb{Thm1.3} {\bf (Liouville theorem)}
Let $M^n$ be a complete noncompact manifold with nonnegative Ricci curvature. Let u be a positive ancient solution, a solution defined in all space and negative time, of  FDE (\ref{Eq3}) for $1-\frac{4}{n+3}<p<1$. If there is a strickly increasing function $L(s)\in C(\R)$  with $L(s)\rightarrow \infty$ as $s\rightarrow \infty$, such that
\begin{eqnarray*}
u(x,t)=o\Big(L(d(x)) +|t|^{1/(1-p)}\Big)
\end{eqnarray*}
near infinity, then $u$ is a constant function on $M^n$.
\et

%\vspace{2mm}

\r
One might see that the growth condition in the spatial direction in Theorem is very weak, since we might choose
$$L(s)=\exp(\exp(\cdots(\exp(s))\cdots))$$
with $l$ $\exp$ for any $l>0$. Note that one might write any positive harmonic function $v(x)$ as a positive global solution of  $\Delta (u^p)=0$ with $u(x)=v(x)^{1/p}$, hence Yau$'$s celebrated Liouville Theorem for positive harmonic functions is a special case of Theorem \ref{Thm1.3} for time dependent positive solutions of FDE, while one couldn't do this for the Heat Equation (see examples in \cite{SZ}).
\er

%\vspace{2mm}

Secondly for PME, i.e. the equation (\ref{Eq3}) with $p>1$, we first obtain the localized Hamilton-type gradient estimate for dimension $n=1$:

%\vspace{2mm}

\t\lb{Thm1.4} Suppose that $u\leq M$ is a positive solution of PME (\ref{Eq3}) in $Q_{R,T} \equiv B(x_0,R)\times [t_0-T , t_0]\subset M^n\times(-\infty,\infty)$  with $n=1$. Let $G(u)=\frac{p}{p-1}u^{p-1}$, $\alpha=\frac{p}{p-1}M^{p-1}(1+\delta)$ with some constant $\delta>0$. Then for any $p>1$, there exists a constant $C(p)$ depending only on p such that
\begin{eqnarray*}
\sup_{(x, t) \in Q_{R/2,T /2}}\frac{|\nabla_x G(u(x,t))|}{\alpha-G(u(x,t))}\leq C(p)\Bigg( \frac{1+\delta}{\delta R}+\frac{1}{\sqrt{M^{p-1}\delta T}}\Bigg).
\end{eqnarray*}

\et

%\vspace{2mm}

An immediate application of the above gradient estimate is the following time-dependent Liouville theorem for PME with $p>1$ on  $\R$,

%\vspace{2mm}

\t\lb{Thm1.5} {\bf (Liouville theorem)}
Let u be a positive ancient solution, a solution defined in all space and negative time, to PME ($p>1$) on  $\R$, such that
\begin{eqnarray}
u(x,t)=o\Big(d(x)^{1/(p-1)} +|t|^{1/(p-1)}\Big)\nn
\end{eqnarray}
near infinity. Then $u$ is a constant.
\et

%\vspace{2mm}

And for PME on a complete Riemannian manifold $(M^n,g)$ with $n\geq 2$, we obtain the localized Hamilton-type gradient estimate:

\t\lb{Thm1.6}  Let $M$ be a complete Riemannian manifold with dimension $n\geq 2$, $Ric(M^n)\geq -k$, $k \geq 0$. Suppose that $u$ is a positive solution of PME (\ref{Eq3}) in $Q_{R,T} \equiv B(x_0,R)\times [t_0-T , t_0]\subset M^n\times(-\infty,\infty)$, with ${\rm range} (u)=[m,M]$. Let $G(u)=\frac{p}{p-1}u^{p-1}$,
$\alpha=\frac{p}{p-1}M^{p-1}(1+\delta)$ with some small constant $0<\delta\leq \frac{4}{n-1}$. If the following {\bf pinch condition} on $m,\; M$ holds
\begin{eqnarray*}
1\leq \Bigg(\frac{M}{m}\Bigg)^{p-1}<\frac{1}{1+\delta}\Big(\frac{4p}{(n-1)(p-1)}+1\Big),
\end{eqnarray*}
then there exists a constant C(n,p) depending only on n and p, and
\begin{eqnarray*}
\gamma=2p-\frac{(n-1)(p-1)}{2}\frac{M^{p-1}(1+\delta)-m^{p-1}}{m^{p-1}}>0,
\end{eqnarray*}
such that
\begin{eqnarray*}
\sup_{(x, t) \in Q_{R/2,T /2}}\frac{|\nabla_x G(u(x,t))|}{\alpha-G(u(x,t))}\leq C(n,p)\Bigg( \frac{\delta +1}{\gamma\delta R} +\frac{1}{\sqrt{\gamma\delta M^{p-1}T}}+\sqrt{\frac{k}{\delta}}  \Bigg).
\end{eqnarray*}
\et

The rest of this paper is organized as follows. In section 2, we first derive a differential inequality for $w=|\nabla G|^2/(\alpha-G)^2$, where we use the technical linear algebra lemma to reduce the use of Cauchy-Schwarz inequality when we derive the differential inequality, while the free use of Cauchy-Schwarz inequality in their proof resulted that Theorem 7 in \cite{MZS} required much stronger conditions. Finally we apply the maximum principle on the differential inequality  to prove Theorem \ref{Thm1.1} following the argument in \cite{SZ} by using the well-known cut-off function by Li and Yau in \cite{LY}. In section 3, we apply the Theorem \ref{Thm1.1} to study the localized Hamilton-type gradient estimates and Liouville-type theorems for FDE and PME and prove Theorem \ref{Thm1.2}-\ref{Thm1.6}. In the appendix, we prove the technical linear algebra lemma.\\

Here and later $Ric(M^n)$ is the Ricci curvature and a manifold is complete if every geodesic extends to infinity. We use the Einstein summation convention for indices $i,\; j,\; k,$ etc. In particular, we use the short-hand notation $g_{ij}^2 = g_{ij} g_{ij} =\sum_{i,j=1}^ng_{ij}^2$. Hereafter we always use $C$ to denote different constants depending on $n$ only.\\

\s{Hamilton-type gradient estimates}

%\vspace{5mm}

%\vspace{2mm}

In this section, we prove Theorem \ref{Thm1.1} and show  that the trick introduced in the fundamental work of Li and Yau \cite{LY} can be adapted to study the Hamilton-type estimate. Our argument can be considered as an improvement of the argument in \cite{MZS}, which can be regarded as a combination of the analysis in \cite{LY} and \cite{SZ}. We define a quantity $w(x, t)$ and start with deriving a differential inequality  on $w(x,t)$, then use the well-known cut-off function of Li-Yau \cite{LY}, to derive the desired bounds.

%\vspace{2mm}

Let $\phi = \ln u$. Since $u$ is a solution to the equation $u_t = \Delta(F(u))$, simple calculation shows
\ee
\phi_t=\Delta(G(u))+\nabla G(u)\cdot\nabla\phi. \lb{Eq2.1}
\eee
Denote $g(\phi) = G(e^{\phi})$, and multiply (\ref{Eq2.1}) by $g'(\phi)$, with some elementary computations, we get
\y
g_t&=&g'\Delta g+|\nabla g|^2,\quad \nabla g=g' \nabla\phi.\\
g'(\phi)&=&G'(e^{\phi})e^{\phi}=F'(u).\\
g''(\phi)&=&F''(e^{\phi})e^{\phi}=F''(u)u.
\ey
Set
\ee
w=w(x,t)=|\nabla\ln(\alpha-g)|^2=\frac{|\nabla g|^2}{(\alpha-g)^2},\lb{Eq2.2}
\eee
and we first derive the differential inequality for $w$, to which we apply the maximum principle.

%\vspace{2mm}

\la\label{La2} $w$ satisfies the following differential inequality:
\y
g'\Delta w-w_t&\geq& L w^2-2g'kw-L_1\nabla g\cdot\nabla w,
\ey
where $L$ and $L_1$ are some functions give in (\ref{Eq2.10}) and (\ref{Eq2.10-1}).
\el

\begin{proof}[\bf Proof of Lemma \ref{La2}:]  After some elementary computations in local orthonormal system as in \cite{MZS}, we get that
\begin{eqnarray}
w_t&=&2\frac{\nabla g\cdot\nabla g_t}{(\alpha-g)^2}+2\frac{|\nabla g|^2g_t}{(\alpha-g)^3}\nn\\
&=& 2\frac{\nabla g\cdot\nabla \big(g'\Delta g+|\nabla g|^2\big)}{(\alpha-g)^2}
+2\frac{|\nabla g|^2\big(g'\Delta g +|\nabla g|^2\big)}{(\alpha-g)^3}\nn\\
&=&\frac{2g'\nabla g\cdot\nabla \Delta g}{(\alpha-g)^2} +\frac{2\Delta g\nabla g\cdot\nabla g'}{(\alpha-g)^2}
+\frac{2\nabla g\cdot\nabla |\nabla g|^2}{(\alpha-g)^2} +\frac{2g'|\nabla g|^2\Delta g}{(\alpha-g)^3}+\frac{2|\nabla g|^4}{(\alpha-g)^3}\nn\\
&=& 2g'\frac{g_jg_{iij}}{(\alpha-g)^2} +2\frac{g''}{g'}\frac{|\nabla g|^2\Delta g}{(\alpha-g)^2}
+4\frac{g_ig_{ij}g_j}{(\alpha-g)^2}+2g'\frac{|\nabla g|^2\Delta g}{(\alpha-g)^3}+2\frac{|\nabla g|^4}{(\alpha-g)^3}\lb{Eq2.3}
\end{eqnarray}

\begin{eqnarray}
w_j=\Big(\frac{g_i^2}{(\alpha-g)^2}\Big)_j =2\frac{g_ig_{ij}}{(\alpha-g)^2}+2\frac{g_i^2g_j}{(\alpha-g)^3}\lb{Eq2.4}
\end{eqnarray}

\begin{eqnarray}
\Delta w&=& w_{jj}=2\Big(\frac{g_ig_{ij}}{(\alpha-g)^2}\Big)_j+2\Big(\frac{g_i^2g_j}{(\alpha-g)^3}\Big)_j\nn\\
&=&2\frac{g_{ij}^2}{(\alpha-g)^2}+2\frac{g_ig_{ijj}}{(\alpha-g)^2}+8\frac{g_ig_{ij}g_j}{(\alpha-g)^3}+2\frac{|\nabla g|^2\Delta g}{(\alpha-g)^3} +6\frac{|\nabla g|^4}{(\alpha-g)^4}\lb{Eq2.5}
\end{eqnarray}
By (\ref{Eq2.3}) and (\ref{Eq2.5}), we obtain that
\begin{eqnarray}
g'\Delta w-w_t&=&2g'\frac{g_{ij}^2}{(\alpha-g)^2}+2g'\frac{g_ig_{ijj}-g_j g_{iij}}{(\alpha-g)^2}+8g'\frac{g_ig_{ij}g_j}{(\alpha-g)^3} \nn\\ &&-4\frac{g_ig_{ij}g_j}{(\alpha-g)^2}+6g'\frac{|\nabla g|^4}{(\alpha-g)^4}-2\frac{g''}{g'}\frac{|\nabla g|^2\Delta g}{(\alpha-g)^2}-2\frac{|\nabla g|^4}{(\alpha-g)^3}\nn
\end{eqnarray}
Bochner¡¯s identity implies that
\begin{eqnarray}
g_ig_{ijj} -g_j g_{iij} = g_j (g_{jii}- g_{iij} ) = R_{ij} g_ig_j=Ric(\nabla g, \nabla g)\nn
\end{eqnarray}
where $R_{ij}$ is the Ricci curvature tensor. Therefore we have
\begin{eqnarray}
g'\Delta w-w_t&=&2g'\frac{g_{ij}^2}{(\alpha-g)^2}+2g'\frac{Ric(\nabla g, \nabla g)}{(\alpha-g)^2} +8g'\frac{g_ig_{ij}g_j}{(\alpha-g)^3}-4\frac{g_ig_{ij}g_j}{(\alpha-g)^2} \nn\\
&&+6g'\frac{|\nabla g|^4}{(\alpha-g)^4}-2\frac{g''}{g'}\frac{|\nabla g|^2\Delta g}{(\alpha-g)^2}-2\frac{|\nabla g|^4}{(\alpha-g)^3}\lb{Eq2.6}
\end{eqnarray}
Recalling (\ref{Eq2.4}), we have
\begin{eqnarray}
\nabla g\cdot\nabla w=2\frac{g_ig_{ij}g_j}{(\alpha-g)^2}+2\frac{|\nabla g|^4}{(\alpha-g)^3}\lb{Eq2.7}
\end{eqnarray}
Adding $\Big(2-\frac{2g'}{\alpha-g}-\eta \frac{g''}{g'}\Big)\times$(\ref{Eq2.7}) with (\ref{Eq2.6}), where $\eta$ is a parameter function to be determined later, we conclude that
\begin{eqnarray}
g'\Delta w-w_t&=&2g'\frac{g_{ij}^2}{(\alpha-g)^2}+2\Big(\frac{2g'}{\alpha-g}-\eta \frac{g''}{g'}\Big)\frac{g_ig_{ij}g_j}{(\alpha-g)^2}-2\frac{g''}{g'}\frac{|\nabla g|^2\Delta g}{(\alpha-g)^2} \nn\\
&&+\Big(2g'+\bigg[2-2\eta \frac{g''}{g'}\bigg](\alpha-g)\Big)\frac{|\nabla g|^4}{(\alpha-g)^4}\nn\\
&&+2g'\frac{Ric(\nabla g, \nabla g)}{(\alpha-g)^2}-\Big(2-\frac{2g'}{\alpha-g}-\eta \frac{g''}{g'}\Big)\nabla g\cdot\nabla w \lb{Eq2.8}
\end{eqnarray}
Denote $f=2g'/(\alpha-g)$, $b=g''/g'$, $A=(g_{ij})$ and $v=\nabla g/|\nabla g|$, then we have $g_{ij}^2=|A|^2$, $g_ig_{ij}g_j=A(v,v)|\nabla g|^2$ and $\Delta g=\tr A$. From definition of $w$, the right side of (\ref{Eq2.8}) can be written as
\begin{eqnarray}
&&2g'\frac{|A|^2}{(\alpha-g)^2}+2(\alpha-g)\Bigg[\Big(f-\eta b\Big)\frac{A(v,v)}{|A|}
-b\frac{\tr A}{|A|}\Bigg]\frac{|A|}{\alpha-g}w \nn\\
&&+(\alpha-g)\Big(f+2-2\eta b\Big)w^2+2g'Ric(v, v)w-\Big(2-f-\eta b\Big)\nabla g\cdot\nabla w \nn\\
&=&2g'\Bigg[\frac{|A|}{(\alpha-g)} +\frac{1}{f}\Bigg(\Big(f-\eta b\Big)\frac{A(v,v)}{|A|}
-b\frac{\tr A}{|A|}\Bigg)w\Bigg]^2\nn\\
&&+(\alpha-g)\Bigg(f+2-2\eta b-\frac{1}{f}\bigg[\Big(f-\eta b\Big)\frac{A(v,v)}{|A|} -b\frac{\tr A}{|A|}\bigg]^2\Bigg)w^2 \nn\\
&&+2g'Ric(v, v)w-\Big(2-f-2\eta b\Big)\nabla g\cdot\nabla w \nn\\
&\geq&\frac{\alpha-g}{f}\Bigg(f^2 +\bigg[2-2\eta b\bigg]f -\bigg[\Big(f-\eta b\Big)\frac{A(v,v)}{|A|}-b\frac{\tr A}{|A|}\bigg]^2\Bigg)w^2 \nn\\
&&+2g'Ric(v, v)w-\Big(2-f-\eta b\Big)\nabla g\cdot\nabla w \lb{Eq2.9}
\end{eqnarray}
Apply Lemma \ref{La1} from Appendix A to (\ref{Eq2.9}), we have
\begin{eqnarray}
g'\Delta w-w_t&\geq&\frac{\alpha-g}{f}\Bigg(f^2 +\Big(2-2\eta b\Big)f -\Big(f-\eta b-b\Big)^2-(n-1)b^2\Bigg)w^2 \nn\\
&&+2g'Ric(v, v)w-\Bigg(2-f-\eta b\Bigg)\nabla g\cdot\nabla w\nn\\
&=&\frac{\alpha-g}{f}\Big(2(1+b)f -\big(n-1+[\eta+1]^2\big)b^2\Big)w^2 \nn\\
&&+2g'Ric(v, v)w-\Bigg(2-f-\eta b\Bigg)\nabla g\cdot\nabla w \nn
\end{eqnarray}
Next we estimate the coefficient $L$ of $w^2$. Pick $\eta=-1$, We can bound $L$ as the following:

\ee
L=(\alpha-g)\Big(2(1+b) -(n-1)\frac{b^2}{f}\Big)\geq (\alpha-g)\gamma >0,\lb{Eq2.10}
\eee
since $2(1+b) -(n-1)\frac{b^2}{f}\geq\gamma>0$ from our condition (C) in Theorem \ref{Thm1.1}. We estimate the coefficient $L_1$ of $\nabla g\cdot\nabla w$ as
\ee
|L_1|= \Big|2-f+b\Big|\leq 2+\tau+f. \lb{Eq2.10-1}
\eee

Hence, from $Ric(M^n)\geq -k$,  we have
\begin{eqnarray}
g'\Delta w-w_t\geq Lw^2-2g'kw-L_1\nabla g\cdot\nabla w \lb{Eq2.12}
\end{eqnarray}
\end{proof}

%\vspace{5mm}

Now we can apply maximum principle to the differential inequality (\ref{Eq2.12}) to prove our gradient estimate (\ref{Eq2}). We will follow \cite{MZS} and \cite{SZ} to use the well-known cut-off function by Li and Yau \cite{LY} to show Theorem \ref{Thm1.1}. We caution the reader that the calculation is not the same as that in \cite{LY} due to the difference of the first-order term.

%\vspace{2mm}

\begin{proof}[\bf Proof of Theorem \ref{Thm1.1}:]
Let $\Psi = \Psi(x, t)$ be a smooth cut-off function supported in $Q_{R,T}$, satisfying the following properties:

\begin{itemize}
\item[]
(1). $\Psi = \Psi(d(x,x_0), t)$, $\Psi=1$ in $Q_{R/2,T/2}$, $0\leq \Psi\leq 1$;

\item[]
(2). $\Psi$ is decreasing as a radial function in the spatial variables;
\item[]
(3). $\frac{\partial_r\Psi}{\Psi^a}\leq \frac{C_a}{R}$, $\frac{\partial_r^2\Psi}{\Psi^a}\leq \frac{C_a}{R^2}$, when $0<a <1$;
\item[]
(4). $\frac{|\partial_t\Psi|^2}{\Psi}\leq \frac{C}{T^2}$
\end{itemize}

Then, from (\ref{Eq2.12}) and a straightforward calculation, one has
\begin{eqnarray}
&&g'\Delta (\Psi w)-(\Psi w)_t+ L_1 \nabla g\cdot \nabla(\Psi w)-2g'\frac{\nabla\Psi}{\Psi}\cdot \nabla(\Psi w)\lb{Eq2.13}\\
&\geq&L\Psi w^2-2kg'\Psi w+g'\Delta\Psi w-\Psi_tw+L_1(\nabla g\cdot \nabla\Psi)w-2g'\frac{|\nabla\Psi|^2}{\Psi}w. \nn
\end{eqnarray}

We obtain the upper bounds for each term of the right-hand side of (\ref{Eq2.13}) as did by Souplet and Zhang \cite{SZ}.
\begin{eqnarray}
\Big|2kg'\Psi w\Big| &\leq& \frac{1}{6}L\Psi w^2+ \frac{6k^2g'^2\Psi}{L} \leq  \frac{1}{6}L\Psi w^2+ Ck^2\frac{g'^2}{L};\lb{Eq2.14}\\
&&\nn\\
\Big|\Psi_t w\Big| &\leq& \frac{1}{6}L\Psi w^2+ \frac{3|\Psi_t|^2}{2L\Psi} \leq  \frac{1}{6}L\Psi w^2+ \frac{C}{T^2}\frac{1}{L}; \lb{Eq2.15}\\
&&\nn\\
\Big|2g'\frac{|\nabla\Psi|^2}{\Psi}w\Big| &\leq& \frac{1}{6}L\Psi w^2+ \frac{6g'^2}{L}\frac{|\nabla\Psi|^4}{\Psi^3}\leq \frac{1}{6}L\Psi w^2+ \frac{C}{R^4}\frac{g'^2}{L};\lb{Eq2.16}\\
&&\nn\\
\Big|L_1(\nabla g\cdot \nabla\Psi)w\Big| &=& \Big|L_1 (\nabla g\cdot \nabla\Psi) \frac{\alpha-g}{|\nabla g|} w^{3/2}\Big|\nn\\
&\leq& \frac{1}{6}L\Psi w^2+ \frac{CL_1^4(\alpha-g)^4}{L^3}\frac{|\nabla\Psi|^4}{\Psi^3}\nn\\
&\leq& \frac{1}{6}L\Psi w^2+ \frac{C}{R^4}\frac{L_1^4(\alpha-g)^4}{L^3};\lb{Eq2.17}
\end{eqnarray}
Here we use Young's inequality,
$$ab\leq \frac{a^p}{p} +\frac{b^q}{q},\qquad \forall p,q>0\; with\; \frac{1}{p}+\frac{1}{q}=1.$$
Furthermore, by the properties of $\Psi$ and the assumption on the Ricci curvature,
one has
\begin{eqnarray}
|-g'\Delta\Psi w|&\leq& \Big|g'\Big[\partial_r^2\Psi+(n-1)\frac{\partial_r\Psi}{r}+\partial_r\Psi\partial_r\ln(\sqrt{\g})\Big]w\Big|\nn\\
&\leq& g'\Big[|\partial_r^2\Psi|+2(n-1)\frac{|\partial_r\Psi|}{R}+k|\partial_r\Psi|\Big]w\nn\\
&\leq& \frac{1}{6}L\Psi w^2+ \frac{9g'^2}{L}\Big[\frac{|\partial_r^2\Psi|^2}{\Psi}+4(n-1)^2\frac{|\partial_r\Psi|^2}{R^2\Psi} +k\frac{|\partial_r\Psi|^2}{\Psi}\Big]\nn\\
&\leq& \frac{1}{6}L\Psi w^2+ C\Big[\frac{1}{R^4}+\frac{1}{R^4} +\frac{k}{R^2}\Big]\frac{g'^2}{L}\lb{Eq2.18}
\end{eqnarray}
Inserting (\ref{Eq2.14})-(\ref{Eq2.18}) into the right-hand side of (\ref{Eq2.13}), we deduce that
\begin{eqnarray}
&&g'\Delta (\Psi w)-(\Psi w)_t+ L_1 \nabla g\cdot \nabla(\Psi w)-2g'\frac{\nabla\Psi}{\Psi}\cdot \nabla(\Psi w)\nn\\
&\geq&\frac{L}{6}\Psi w^2- C\Bigg[\frac{L_1^4(\alpha-g)^4}{L^3}\frac{1}{R^4}+\frac{g'^2}{L}\frac{1}{R^4}+\frac{1}{L}\frac{1}{T^2}+\frac{g'^2}{L}k^2 \Bigg] \nn\\
&\geq&\frac{L}{6}\Bigg[\Psi w^2- C\Bigg[\frac{L_1^4(\alpha-g)^4}{L^4}\frac{1}{R^4}+\frac{g'^2}{L^2}\frac{1}{R^4}+\frac{1}{L^2}\frac{1}{T^2}+\frac{g'^2}{L^2}k^2 \Bigg]\Bigg]\lb{Eq2.19}
\end{eqnarray}

Recalling that $L\geq 2\gamma(\alpha-g)>0$, $L_1\leq 2+\tau+f$ and $f=g'/(\alpha-g)$, therefore,
\begin{eqnarray}
&&g'\Delta (\Psi w)-(\Psi w)_t+ L_1 \nabla g\cdot \nabla(\Psi w)-2g'\frac{\nabla\Psi}{\Psi}\cdot \nabla(\Psi w)\nn\\
&\geq&\frac{L}{6}\Bigg[\Psi w^2- C(\gamma)\Bigg( \frac{(1+f+\tau )^4}{R^4}
+\frac{1}{(\alpha-g)^2}\Big[\frac{g'^2}{R^4}+\frac{1}{T^2}+g'^2 k^2\Big] \Bigg)\Bigg]\nn\\
&\geq&\frac{L}{6}\Bigg[\Psi w^2- C(\gamma)\Bigg( \frac{(1+f+\tau )^4}{R^4}
+\frac{1}{(\alpha-g)^2}\frac{1}{T^2}+f^2 k^2\Bigg)\Bigg]. \lb{Eq2.20}
\end{eqnarray}

Suppose that the maximum of $(\Psi w)$ is reached at $(x_1, t_1)\in Q_{R,T}$. By \cite{LY}, we can assume,
without loss of generality, that $x_1$ is not in the cut-locus of $M$. Then at this point,
one has $\Delta (\Psi w)\leq 0$, $(\Psi w)_t \geq 0$ and $\nabla(\Psi w) = 0$. Recalling that $\alpha-g\geq \delta$, $0<g'\leq K$ and $f\leq K/\delta$, we have
\begin{eqnarray}
\Big(\Psi w^2\Big)(x_1,t_1) \leq C(\gamma)\Bigg( \frac{(\delta+K+\tau\delta )^4}{\delta^4R^4}+\frac{1}{\delta^2}\frac{1}{T^2} +\Big(\frac{K}{\delta}\Big)^2 k^2 \Bigg)\nn
\end{eqnarray}

By assumption, the maximum of $(\Psi w)$ is reached at $(x_1, t_1)\in Q_{R,T}$, which implies that for any $(x, t)\in Q_{R,T}$
\begin{eqnarray}
\Big(\Psi w\Big)^2(x,t)&\leq&\Big(\Psi w\Big)^2(x_1,t_1)\leq\Big(\Psi w^2\Big)(x_1,t_1)\nn\\
&\leq& C(\gamma)\Bigg(\frac{(\delta+K+\tau\delta )^4}{\delta^4R^4} +\frac{1}{\delta^2T^2}+\Big(\frac{K}{\delta}\Big)^2 k^2 \Bigg) \nn
\end{eqnarray}
Noticing that $\Psi(x, t) = 1$ in $Q_{R/2,T /2}$ and $w = |\nabla g|^2/(\alpha-g)^2$, we finally have proven
\begin{eqnarray}
\frac{|\nabla g|^2}{(\alpha-g)^2}\leq C(\gamma)\Bigg(\frac{(\delta+K+\tau \delta)^2}{\delta^2R^2} +\frac{1}{\delta T }+\frac{K}{\delta} k \Bigg) \nn
\end{eqnarray}
which is exactly what the conclusion of Theorem \ref{Thm1.1} is.\\
\end{proof}

\s{Gradient Estimates and Liouville Theorems for FDE and PME}

%\vspace{2mm}

In this section, we study the heat equation, FDE and PME on a complete Riemannian manifold, and derive the localized Hamilton-type gradient estimates by applying our Theorem \ref{Thm1.1}, then prove some time-dependent Liouville Theorems for FDE and PME on noncompact complete manifolds with nonnegative Ricci curvature.  As  a corollary, we obtain Yau$'$s celebrated Liouville theorem for positive harmonic functions: {\it any positive harmonic function on a noncompact manifold with nonnegative Ricci curvature is a constant function.}\\

\ss{Heat Equations:}

Let $M^n$ be a complete Riemannian manifold with dimension $n\geq 1$, $Ric(M^n)\geq -k$, $k \geq 0$. Suppose that $u\leq M$ is a positive solution of the heat equation
$$ u_t=\Delta u $$
in $Q_{R,T} \equiv B(x_0,R)\times [t_0-T , t_0]\subset M^n\times(-\infty,\infty)$. We may look as the case of $F(u)=u$ in Theorem \ref{Thm1.1}. Choose $G(s)=\ln s$ and let $\alpha=1+\ln M$, $K=1$, $\tau=1$ and $\gamma=2$ in Theorem \ref{Thm1.1}, we obtain the main result of  Souplet and  Zhang  \cite{SZ}, {\bf Theorem C}, as a direct corollary of our Theorem  \ref{Thm1.1}. \\

Using the localized Hamilton-type gradient estimates for positive solutions of the heat equation,  Souplet and  Zhang  \cite{SZ})  proved the following time-dependent Liouville Theorem:

\c {\bf (Theorem 1.2 in \cite{SZ})} Let $M$ be a complete, noncompact manifold with nonnegative Ricci curvature. Then the following conclusions hold.

(a) Let $u$ be a positive ancient solution to the heat equation (that is, a solution defined in all space and negative time) such that
$u(x, t) = e^{o(d(x)+\sqrt{|t|})}$ near infinity. Then u is a constant.

(b) Let $u$ be an ancient solution to the heat equation such that $u(x, t) = o([d(x)+\sqrt{|t|}])$ near infinity. Then u is a constant.
\ec

As discussed in \cite{SZ}, one Could not expect that Yau¡¯s Liouville theorem would still hold for positive ancient or eternal solutions to the heat equation. Both growth conditions of the above theorem in the spatial direction are sharp, by some simple examples, i.e., for (a), let $u(x,t)=e^{x+t}$ on $x\in \R$, and for (b), let $u(x,t)=x$ on $x\in \R$. Hence one couldn't obtain  Yau$'$s celebrated Liouville theorem for positive harmonic functions directly from the above time-dependent Liouville Theorem for heat equation.\\

%\vspace{2mm}

\ss{Fast Diffusion Equation:}

Let $M^n$ be a complete Riemannian manifold with dimension $n\geq 1$, $Ric(M^n)\geq -k$, $k \geq 0$. Here we consider the Fast Diffusion Equation (FDE for short)
$$ u_t=\Delta (u^p),\qquad \qquad 0<p<1 $$
on $M^n\times(-\infty,\infty)$. We have the following localized Hamilton-type gradient estimates for the positive solution on $M^n\times(-\infty,\infty)$:

\t\lb{Thm3.1} Suppose that $u\leq M$ is a positive solution of the Fast Diffusion Equation
$$ u_t=\Delta (u^p) $$
in $Q_{R,T} \equiv B(x_0,R)\times [t_0-T , t_0]\subset M^n\times(-\infty,\infty)$, where
\begin{eqnarray}
1-\frac{4}{n+3}< p< 1,  & for\; n\geq 1.\lb{PMEq<1}
\end{eqnarray}

Then there exists a constant C depending only on n and p such that
\begin{eqnarray*}
\sup_{(x, t)\in Q_{R/2,T /2}}\frac{|\nabla u(x,t)|}{u(x,t)}\leq  C\Big(\frac{1}{R}+\frac{M^{(1-p)/2}}{\sqrt{T}}+\sqrt{k} \Big).
\end{eqnarray*}
When $n=1$, the Ricci curvature lower bound $k$ vanishes.

\et

%\vspace{2mm}

\begin{proof}[\bf Proof of Theorem \ref{Thm3.1}:] For $0 < p <1$, we have $F(s) = s^p$, and choose $G(s)=\frac{p}{p-1}s^{p-1}$, and  let $\alpha = 0$ and $\tau \geq 1-p$ in Theorem \ref{Thm1.1}. We have $F'(s)=ps^{p-1}$, $F''(s)=p(p-1)s^{p-2}$, then condition (C) in Theorem \ref{Thm1.1} becomes
\begin{eqnarray*}
4p-(n-1)(1-p)>0,
\end{eqnarray*}
which is equivalent to,
\begin{eqnarray*}
1-\frac{4}{n+3} < p <1.
\end{eqnarray*}
Let $\gamma=\big((n+3)p-(n-1)\big)/2$ for given $p$ in admission range (\ref{PMEq<1}). Follow the proof of Theorem \ref{Thm1.1}, we have
\begin{eqnarray}
\left\{
\begin{array}{rl}
& g=\frac{p}{p-1}u^{p-1},\; g'=pu^{p-1},\; g''=p(p-1)u^{p-1}\\
&\\
& f=\frac{2g'}{-g}=2(1-p), \quad b=\frac{g''}{g'}=p-1, \\
&\\
& L=\gamma \frac{p}{1-p}u^{p-1}, \quad L_1=3-p.
\end{array} \right.\nn
\end{eqnarray}
We follow the proof of Theorem \ref{Thm1.1} until (\ref{Eq2.19}) with the modification that the constant $C(n,p)$ here depends only on $n$ and $p$.  Inserting the above quantities into (\ref{Eq2.19}), we have
\begin{eqnarray}
&&2g'\Delta (\Psi w)-(\Psi w)_t+ L_1 \nabla g\cdot \nabla(\Psi w)-2g'\frac{\nabla\Psi}{\Psi}\cdot \nabla(\Psi w)\nn\\
&\geq&\frac{L}{6}\Bigg[\Psi w^2- C(n,p)\Bigg( \frac{1}{R^4}
+\frac{1}{u^{2(p-1)}}\frac{1}{T^2}+ k^2\Bigg)\Bigg]. \nn
\end{eqnarray}
By the same maximum argument in proof of Theorem \ref{Thm1.1}, we have
\y
\Psi w^2\leq C(n,p)\Bigg( \frac{1}{R^4}+\frac{M^{2(1-p)}}{T^2}+ k^2\Bigg),
\ey
which implies
\begin{eqnarray*}
\frac{|\nabla g|^2}{(-g)^2}\leq C(n,p)\Bigg( \frac{1}{R^2}+\frac{M^{1-p}}{T}+ k\Bigg).
\end{eqnarray*}
Then the conclusion of Theorem \ref{Thm3.1} follows easily from the fact that
\begin{eqnarray*}
\frac{|\nabla g|}{-g}=\frac{|\nabla G(u)|}{-G(u)}=(1-p)\frac{|\nabla u|}{u}\nn
\end{eqnarray*}
\end{proof}

%\vspace{2mm}

An immediate application of the above gradient estimates is the following time-dependent Liouville theorem for FDE on a complete noncompact manifold with nonnegative Ricci curvature:

%\vspace{2mm}

\t\lb{Thm3.2} {\bf (Liouville theorem)}
Let $M^n$ be a complete, noncompact manifold with nonnegative Ricci curvature. Let u be a positive ancient solution, a solution defined in all space and negative time, of the Fast Diffusion Equation for $1-\frac{4}{n+3}<p<1$ , and $L(s)\in C(\R)$ be any strickly increasing function with $L(s)\rightarrow \infty$ as $s\rightarrow \infty$, such that
\begin{eqnarray}
u(x,t)=o\Big(L(d(x)) +|t|^{1/(1-p)}\Big)\nn
\end{eqnarray}
near infinity. Then $u$ is a constant.
\et

%\vspace{2mm}

\begin{proof}[\bf Proof of Theorem \ref{Thm3.2}:]  Since $L(s)$ is a strickly increasing function with $L(s)\rightarrow \infty$ as $s\rightarrow \infty$, there is an inverse function $H(s)$ of $L(s)$ which is also a strickly increasing function with $H(s)\rightarrow \infty$ as $s\rightarrow \infty$. Fixing $(x_0, t_0)$ in space-time and using Theorem \ref{Thm3.1} for $u$ on the cube $Q(\frac{1}{2}H(R^{2/(1-p)}), R^2)=B(x_0,\frac{1}{2}H(R^{2/(1-p)}))\times [t_0-R^2, t_0]$, and the $M$ in Theorem 3.1 is the maximum value on the double cube $Q(H(R^{2/(1-p)}), 2R^2)$, by  our assumption on the growth condition of the function $u$ at infinity,
\ee
M_{H(R^{2/(1-p)}),2R^2}=o\Big(L(H(R^{2/(1-p)}))+R^{2/(1-p)}\Big)=o(R^{2/(1-p)}).\nn
\eee
Hence by Theorem \ref{Thm3.1}, we have
\ee
\frac{|\nabla u(x_0,t_0)|}{u(x_0,t_0)}\leq  C\Big(\frac{1}{\frac{1}{2}H(R^{2/(1-p)})}+\frac{1}{R}o(R)\Big), \nn
\eee
Letting $R \rightarrow\infty$, it follows that $|\nabla u(x_0, t_0)| = 0$. Since $(x_0, t_0)\in M^n\times \R$ is arbitrary, one
sees that $u$ must be a constant function.
\end{proof}

As a corollary of the above time dependent Liouville theorem for FDE, we obtain Yau$'$s celebrated Liouville theorem for positive harmonic functions on a complete, noncompact manifold with nonnegative Ricci curvature:

\c\label{Cor3.3}{\bf (Yau$'$s Liouville theorem for positive harmonic functions)}
Any positive harmonic function on a noncompact manifold with nonnegative Ricci curvature is a constant function.
\ec

\begin{proof}[\bf Proof of Corollary \ref{Cor3.3}:]  The proof  follows  immediately from the above time dependent Liouville theorem. Let $v$ be a positive harmonic function. Choose a $p$ with $1-\frac{4}{n+3}<p<1$, then $u(x)=v(x)^{1/p}$ is a positive solution of $\Delta (u^p)=0$, which can be considered as a time independent positive solution of the corresponding Fast Diffusion Equation. Define
$$L(s)=s\max_{d(x)\leq s} u(x) +s$$
It is easy to see that $L(s)$ is a strickly increasing function with $L(s)\rightarrow \infty$ as $s\rightarrow \infty$, and $u(x,t)=o(L(d(x)))$ near infinity. Follow from the above Theorem, $u$ must be a constant function, which implies $v$ must be a constant function.
\end{proof}

\ss{Porous Media Equations:}

%\vspace{2mm}

Let $M^n$ be a complete Riemannian manifold with dimension $n\geq 1$, $Ric(M^n)\geq -k$, $k \geq 0$. Here we consider Porous Media Equation (PME for short)
$$ u_t=\Delta (u^p),\qquad \qquad p>1 $$
on $M^n\times(-\infty,\infty)$.\\

For dimension $n=1$, we have the following localized Hamilton-type gradient estimates for the positive solution on $M^n\times(-\infty,\infty)$:

\t\lb{Thm3.3} Suppose that $u\leq M$ is a positive solution of PME (\ref{Eq3}) in $Q_{R,T} \equiv B(x_0,R)\times [t_0-T , t_0]\subset M^n\times(-\infty,\infty)$  with $n=1$. Let $G(u)=\frac{p}{p-1}u^{p-1}$, $\alpha=\frac{p}{p-1}M^{p-1}(1+\delta)$ with some constant $\delta>0$. Then for any $p>1$, there exists a constant $C(p)$ depending only on p such that
\begin{eqnarray*}
\sup_{(x, t) \in Q_{R/2,T /2}}\frac{|\nabla_x G(u(x,t))|}{\alpha-G(u(x,t))}\leq C(p)\Bigg( \frac{1+\delta}{\delta R}+\frac{1}{\sqrt{M^{p-1}\delta T}}\Bigg).
\end{eqnarray*}

\et

%\vspace{2mm}

\begin{proof}[\bf Proof of Theorem \ref{Thm3.3}:]  We have $F(s) = s^p$ and choose $G(s)=\frac{p}{p-1}s^{p-1}$ in Theorem \ref{Thm1.1}. We have $F'(s)=ps^{p-1}$, $F''(s)=p(p-1)s^{p-2}$,and   $K=pM^{p-1}$, $\tau \geq p-1$, and $\alpha$ as defined in above, then condition (C) in Theorem 1.1 becomes
\begin{eqnarray*}
 2p \geq \gamma>0.
\end{eqnarray*}
Condition (C) is satisfied if we choose $\gamma=2p$. Similar to the proof of Theorem \ref{Thm3.1}, we have
\begin{eqnarray}
\left\{
\begin{array}{rl}
& g=\frac{p}{p-1}u^{p-1},\; g'=pu^{p-1},\; g''=p(p-1)u^{p-1}\\
&\\
& f=\frac{2g'}{\alpha-g}=2(p-1)\frac{u^{p-1}}{M^{p-1}(1+\delta)-u^{p-1}}\leq\frac{2(p-1)}{\delta}\\
&\\
& b=\frac{g''}{g'}=p-1\\
&\\
& L=\gamma(\alpha-g)=\frac{2p^2}{1-p}\Big[M^{p-1}(1+\delta)-u^{p-1}\Big]\geq \frac{2p^2}{1-p}M^{p-1}\delta\\
&\\
& L_1=1+p+f\leq 1+p+2(p-1)\frac{u^{p-1}}{M^{p-1}(1+\delta)-u^{p-1}}\leq 1+p+\frac{2(p-1)}{\delta}.
\end{array} \right.\nn
\end{eqnarray}
By the same argument in the proof of Theorem  \ref{Thm3.1}, we have
\y
\Psi w^2\leq C(p)\Bigg( \frac{L_1^4}{R^4}+\frac{1}{L^2T^2}\Bigg)
\leq C(p)\Bigg( \frac{(1+\frac{1}{\delta})^4}{R^4}+\frac{1}{(M^{p-1}\delta)^2 T^2}\Bigg),
\ey
which implies
\begin{eqnarray*}
\frac{|\nabla g|^2}{(\alpha-g)^2}\leq C(p)\Bigg( \frac{(1+\frac{1}{\delta} )^2}{ R^2}+\frac{1}{M^{p-1}\delta T}\Bigg).
\end{eqnarray*}
Then the conclusion of Theorem \ref{Thm3.3} follows easily from the fact that
\begin{eqnarray*}
\frac{|\nabla g|}{\alpha-g}=\frac{|\nabla G(u)|}{\alpha-G(u)}.
\end{eqnarray*}
\end{proof}

%\vspace{2mm}

An immediate application of this theorem is the following time-dependent Liouville theorem for PME on  $\R$.

%\vspace{2mm}

\t\lb{Thm3.4} {\bf (Liouville theorem)}
Let u be a positive ancient solution, a solution defined in all space and negative time, to the Porous Medium Equation ($p>1$) on $\R$, such that
\begin{eqnarray}
u(x,t)=o\Big(d(x)^{1/(p-1)} +|t|^{1/(p-1)}\Big)\nn
\end{eqnarray}
near infinity. Then $u$ is a constant.
\et

%\vspace{2mm}

\begin{proof}[\bf Proof of Theorem \ref{Thm3.4}:]  By our assumption, the function $u$ satisfies $u(x,t) = o(d(x)^{1/(p-1)} + |t|^{1/(p-1)})$ near infinity. Fixing $(x_0, t_0)$ in space-time and using Theorem \ref{Thm3.3} for $u$ on the cube $Q(R, R)=B(x_0,R)\times [t_0-R, t_0]$, we have
\begin{eqnarray*}
\frac{|\nabla u(x_0,t_0)|}{u(x_0, t_0)}\leq C(p)\Big(\frac{1+\delta}{\delta R}+\frac{1}{\sqrt{M^{p-1}\delta T}}\Big)M^{p-1}\leq C(p,\delta)\Big( \frac{M^{p-1}}{R}+\sqrt{\frac{M^{p-1}}{T}}\Big),
\end{eqnarray*}
where $M$ is the maximum value on the double cube $Q(2R, 2R)$, by  our assumption on the growth condition of the function $u$ at infinity,
\ee
M_{2R,2R}=o\Big(R^{1/(p-1)}+R^{1/(p-1)}\Big)=o\big(R^{1/(p-1)}\big).\nn
\eee
Hence we have
\ee
\frac{|\nabla u(x_0,t_0)|}{u(x_0,t_0)}\leq C(p,\delta)\Bigg( \frac{M^{p-1}}{R}+\sqrt{\frac{M^{p-1}}{R}}\Bigg)=o(1), \nn
\eee
Letting $R \rightarrow\infty$, it follows that $|\nabla u(x_0, t_0)| = 0$. Since $(x_0, t_0)$ is arbitrary, one
sees that $u $ must be a constant function.
\end{proof}

%\vspace{2mm}

For dimension $n\geq 2$, we have the following localized Hamilton-type gradient estimates for the positive solution of PME on $M^n\times(-\infty,\infty)$:

\t\lb{Thm3.5}  Let $M$ be a complete Riemannian manifold with dimension $n\geq 2$, $Ric(M^n)\geq -k$, $k \geq 0$. Suppose that $u$ is a positive solution of the Porous Media Equation
$$ u_t=\Delta (u^p) $$
in $Q_{R,T} \equiv B(x_0,R)\times [t_0-T , t_0]\subset M^n\times(-\infty,\infty)$, with ${\rm range} (u)=[m,M]$. Let $G(u)=\frac{p}{p-1}u^{p-1}$,
$\alpha=\frac{p}{p-1}M^{p-1}(1+\delta)$ with some small constant $0<\delta\leq \frac{4}{n-1}$. If the following {\bf pinch condition} on $m,\; M$ holds
\begin{eqnarray}
1\leq \Bigg(\frac{M}{m}\Bigg)^{p-1}<\frac{1}{1+\delta}\Big(\frac{4p}{(n-1)(p-1)}+1\Big), \lb{Eq3.2}
\end{eqnarray}
then there exists a constant C(n,p) depending only on n and p, and
\begin{eqnarray}
\gamma=2p-\frac{(n-1)(p-1)}{2}\frac{M^{p-1}(1+\delta)-m^{p-1}}{m^{p-1}}>0,\nn
\end{eqnarray}
such that
\begin{eqnarray}
\frac{|\nabla_x G(u(x,t))|}{\alpha-G(u(x,t))}\leq C(n,p)\Bigg( \frac{\delta +1}{\gamma\delta R} +\frac{1}{\sqrt{\gamma\delta M^{p-1}T}}+\sqrt{\frac{k}{\delta}}  \Bigg). \nn
\end{eqnarray}
for all  $(x, t)$ in $Q_{R/2,T /2}$.

\et

%\vspace{2mm}

\begin{proof}[\bf Proof of Theorem \ref{Thm3.5}:]  We have $F(s) = s^p$ and choose $G(s)=\frac{p}{p-1}s^{p-1}$ in Theorem \ref{Thm1.1}. We have $F'(s)=ps^{p-1}$, $F''(s)=p(p-1)s^{p-2}$,and   $K=pM^{p-1}$, and $\alpha,$ as defined in above, then Condition (A) and (B) are satisfied and Condition (C) in Theorem 1.1 becomes
\begin{eqnarray*}
\left\{
\begin{array}{rl}
&\tau \geq p-1, \\
& \\
& 2p-\frac{(n-1)(p-1)}{2}\frac{M^{p-1}(1+\delta)-s^{p-1}}{s^{p-1}} \geq \gamma>0, \quad \forall\; s\in [m,M]
\end{array} \right.
\end{eqnarray*}

which is equivalent to,
\begin{eqnarray*}
\left\{
\begin{array}{rl}
&\tau \geq p-1, \\
 & \\
&2p-\frac{(n-1)(p-1)}{2}\frac{M^{p-1}(1+\delta)-m^{p-1}}{m^{p-1}} >0
\end{array} \right.
\end{eqnarray*}
Let $\tau =p-1$ and $\gamma=2p-\frac{(n-1)(p-1)}{2}\frac{M^{p-1}(1+\delta)-m^{p-1}}{m^{p-1}}>0$, the above condition as $0<\delta\leq \frac{4}{n-1}$ is equivalent to our pinch condition (\ref{Eq3.2}).

Similar to the proof of Theorem \ref{Thm3.3}, we have
\begin{eqnarray}
\left\{
\begin{array}{rl}
& g=\frac{p}{p-1}u^{p-1},\; g'=pu^{p-1},\; g''=p(p-1)u^{p-1}\\
&\\
& f=\frac{2g'}{\alpha-g}=2(p-1)\frac{u^{p-1}}{M^{p-1}(1+\delta)-u^{p-1}}\leq\frac{2(p-1)}{\delta}\\
&\\
& b=\frac{g''}{g'}=p-1\\
&\\
& L=\gamma(\alpha-g)=\frac{p\gamma}{1-p}\Big[M^{p-1}(1+\delta)-u^{p-1}\Big]\geq \frac{p\gamma}{1-p}M^{p-1}\delta\\
&\\
& L_1=1+p+f\leq 1+p+2(p-1)\frac{u^{p-1}}{M^{p-1}(1+\delta)-u^{p-1}}\leq 1+p+\frac{2(p-1)}{\delta}.
\end{array} \right.\nn
\end{eqnarray}

By the same argument in the proof of Theorem \ref{Thm3.3}, we have
\begin{eqnarray}
\Psi w^2\leq C(n,p)\Bigg( \frac{(\delta+1)^4}{\gamma^4\delta^4R^4}+\frac{1}{\gamma^2\delta^2M^{2(p-1)}T^2}+\frac{k^2}{\delta^2}\Bigg)\nn
\end{eqnarray}
which implies
\begin{eqnarray}
\frac{|\nabla g|^2}{(\alpha-g)^2}\leq C(n,p)\Bigg( \frac{(\delta+1)^2}{\gamma^2\delta^2R^2}+\frac{1}{\gamma \delta M^{p-1}T}+\frac{k}{\delta}\Bigg) \nn
\end{eqnarray}
Then the conclusion of Theorem \ref{Thm3.5} follows easily from the fact that
\begin{eqnarray}
\frac{|\nabla g|}{\alpha-g}=\frac{|\nabla G(u)|}{\alpha-G(u)}\nn
\end{eqnarray}
\end{proof}

%\vspace{2mm}

{\bf Ackowledgement}

\vspace{2mm}

The work was initiated when the author visited McGill University and Centre de recherches mathematiques in Montreal as a CRM Postdoctoral Fellow for {\it Thematic Semester on Dynamical Systems and Evolution Equations}, January-June 2008. The author would like to thank professor Pengfei Guan, Dimitry Jakobson  and John Toth for help during his visit at McGill and CRM, and would also like  to thank Dr. Junfang Li for many useful discussions and comments on this paper.\\

%\vspace{2mm}

%% The Appendices part is started with the command \appendix;
%% appendix sections are then done as normal sections

\appendix{\bf Appendix A}
\vspace{2mm}

Here we show a technical linear algebra lemma for symmetric matrix.

%\vspace{2mm}

\begin{lemma}\label{La1} Let $A=(a_{ij})$ be a nonzero $n\times n$ symmetric matrix with eigenvalues $\{\lambda_k\}$, for any $a,\; b\in\R$, one has the following properties:
\begin{eqnarray}
&(a).& |A|^2=\sum_{i,j=1}^n a_{ij}^2=\tr (A A^T)=\sum_{k=1}^n \lambda_k^2;\nn\\
&& \nn \\
&(b).& \max_{|v|=1}\big(aA+b\tr A I_n\big)(v,v)=a\lambda_i+b\sum_{k=1}^n \lambda_k, \; for \; some \; 1\leq i\leq n \nn\\
&& \min_{|v|=1}\big(aA+b\tr A I_n\big)(v,v)=a\lambda_j+b\sum_{k=1}^n \lambda_k, \;  for \; some \; 1\leq j\leq n.\nn\\
&& \nn \\
&(c).& \max_{A\in \Sa(n);|v|=1}\Big[\frac{aA+b\tr A I_n}{|A|}(v,v) \Big]^2=(a+b)^2+(n-1)b^2.\nn
\end{eqnarray}
\end{lemma}

%\vspace{2mm}

\begin{proof}[\bf Proof of Lemma \ref{La1}:] (a) follows from direct computation and $A$ symmetry.

%\vspace{2mm}

(b) follows from the facts that
$$\max_{|v|=1}\big(aA+b\tr A I_n\big)(v,v), \qquad \min_{|v|=1}\big(aA+b\tr A I_n\big)(v,v)$$
are the maximal and minimal eigenvalue of $aA+b\tr A I_n$, and the eigenvalues of $aA+b\tr A I_n$ are $\{a\lambda_i+b\sum_{k=1}^n \lambda_k\}_{i=1}^n$.

%\vspace{2mm}

To prove (c), apply (a) and (b), we have
\y
\max_{A\in \Sa(n);|v|=1}\Big[\frac{aA+b\tr A I_n}{|A|}(v,v) \Big]^2=\max_{\sum_{k=1}^n \lambda_k^2=1}[a\lambda_1+b\sum_{k=1}^n \lambda_k]^2.
\ey
By Lagrangian multiplier method in Calculus, the extrema of
$$f(\lambda_1,\cdots,\lambda_n)=a\lambda_1+b\sum_{k=1}^n \lambda_k$$
under constrain $\sum_{k=1}^n \lambda_k^2=1$ are
\begin{eqnarray*}
-[(a+b)^2+(n-1)b^2]^{1/2} \leq a\lambda_1+b\sum_{k=1}^n \lambda_k\leq [(a+b)^2+(n-1)b^2]^{1/2}.
\end{eqnarray*}
\end{proof}

%% \section{}
%% \label{}

%% References
%%
%% Following citation commands can be used in the body text:
%% Usage of \cite is as follows:
%%   \cite{key}         ==>>  [#]
%%   \cite[chap. 2]{key} ==>> [#, chap. 2]
%%

%% References with bibTeX database:

\bibliographystyle{elsarticle-num}

%\bibliography{<your-bib-database>}

%% Authors are advised to submit their bibtex database files. They are
%% requested to list a bibtex style file in the manuscript if they do
%% not want to use elsarticle-num.bst.

%% References without bibTeX database:

%\begin{thebibliography}{00}

%% \bibitem must have the following form:
%%   \bibitem{key}...
%%

% \bibitem{}

\end{document}